\documentclass{bull-l}

\theoremstyle{definition}

\theoremstyle{remark}

\numberwithin{equation}{section}

\begin{document}
\boldmath
\title{Nonlocal Mathematics}
% Remove or comment out any unused author tags.
% author one information
\author{Mukul Patel}
\address{Southwestern College, Dept. of Mathematics\\16389 Saratoga, San
Leandro CA 94578}
%\curraddr{}
\email{math@mailroom.com}

%\thanks{Thanks to all the sponsors of this work.}

% author two information
%\author{}
%\address{}
%\curraddr{}
%\email{}
%\thanks{}

\subjclass{Primary 00, Secondary 11, 14, 16, 19, 32, 35, 53, 55, 57R, 58, 81,
83}
\date{February 25, 2000}
% at present the "communicated by" line appears only in ERA and PROC
%\commby{}

\dedicatory{Dedicated to Prof. Phillip E. Parker}
\begin{abstract}
This is a research report on what is best termed `nonlocal' methods in
mathematics. (This is not to be confused with global analysis.) The nonlocal
formulation of  physics in \cite{principia} points to a fresh viewpoint in
mathematics: the \textbf{nonlocal} viewpoint. It involves analyzing objects of
geometry and analysis
using nonlocal methods, as opposed to the classical local methods, e.g.,
Newton's calculus. It also involves analyzing new nonlocal geometries and
nonlocal
analytical objects, i.e. nonlocal fields. In geometry, we introduce and study
(nonlocal) forms, differentials, integrals, connections, curvatures, holonomy,
$G$-structures, etc. In analysis, we analyze local fields using nonlocal
methods (semilocal analysis); nonlocal fields nonlocally (of course); and the
connection between nonlocal
linear analysis and local \emph{nonlinear} analysis. Analysis and geometry are
next
synthesized to yield nonlocal (hence noncommutative) homology, cohomology,
de Rham theory, Hodge theory, Chern-Weil theory, $K$-theory (called
$N$-theory)
and index theory. Applications include theorems such as nonlocal-
noncommutative Riemann-Roch,
Gauss-Bonnet, Hirzebruch signature, etc.
The nonlocal viewpoint is also investigated in algebraic geometry, analytic
geometry, and, to some extent, in arithmetic geometry, resulting in powerful  
bridges  between the classical and nonlocal-noncommutative
aspects.

\end{abstract}
\maketitle
\section*{Introduction}
It is propounded in \cite{principia} that the basic fields of physics are
nonlocal (Sec.\ref{geometry}) and their local aspects are nonlocally related.
Thus,
laws of physics are formulated using a nonlocal calculus. This circumstance
has led to an elaborate and fruitful development of what we term
\textbf{nonlocal}
geometry and analysis. These are not to be confused with global analysis and
global aspects of differential geometry. Our fields are nonlocal in the sense
that they are
not defined point-wise, but rather at ordered pairs of points and their values
are homomorphisms from the tangent space of the first point to that of the
second. This initiates three lines of thought. First one leads to a nonlocal
differential geometry (Sec.\ref{geometry}); the second one leads to a nonlocal
analysis of usual, local (point-wise), fields (Sec.\ref{analysis}); and the
third leads to an analysis of the nonlocal
fields themselves (Sec.\ref{analysis}).

Although all these considerations arose out of the problems of physics,
the work \cite{nang}, the subject matter of the present report, naturally
focuses on the  mathematics itself. The threads mentioned above (\cite{nang},
Parts I \& II) come together
and interact fruitfully, to exhibit a very exciting interplay. The outcome
(\cite{nang}, Part III, and Sec.\ref{synthesis} below) is
a set of nonlocal theories: noncommutative cohomology,
nonlocal de Rham theory, nonlocal Hodge theory, nonlocal Chern-Weil theory,
nonlocal (noncommutative) $K$-theory, and nonlocal index theory. This last
yields nonlocal analogs of various classical theorems such as 
Gauss-Bonnet-Chern,
Hirzebruch-Riemann-Roch, Hirzebruch signature, Lefschetz fixed point, etc.
Nonlocal methods carry over also to algebraic geometry and
analytic geometry: deep connections are developed between commutative
algebraic
geometry and the geometry of non-commutative rings (\cite{nang}, Part IV, and
Sec.\ref{ag} below). Analytic geometry, too, benefits from nonlocal methods
and accompanying noncommutative
cohomology of sheaves (\cite{nang}, Part V, and Sec.\ref{ang} below), and so
does arithmetic geometry (\cite{nang}, Part VI, and Sec.\ref{nt} below).

\section{Geometry}\label{geometry}
A \textbf {nonlocal vector (resp. affine) field} $\omega$ on a manifold $X$ is
an assignment to each ordered pair of points $(x,y)$ a linear (resp. affine)
homomorphism $\omega_{xy}$ from the tangent space $T_x$ to the tangent space
$T_y.$ We call such a field a (nonlocal) $1$-form. Similarly, we can define a
\textbf{(nonlocal) $p$-form} to be an assignment of a set of homomorphisms 
along the
edges of
each singular $p$-cube in $X.$ 

A $1$-form with invertible homomorphisms can also be thought of as a
\textbf{nonlocal connection:} a connection, because it provides a means of
comparing
vectors at
different points; \textbf{nonlocal} in the sense that the comparison is of
vectors at \emph{distinct} points, as opposed to the usual connection from
differential geometry---which is
essentially a means of comparing vectors at a point with those in its
\emph{infinitesimal} vicinity. The latter, too, enables comparison of vectors at
distinct points, but
the comparison depends on the path along which vectors are parallelly
transported---which again underscores the local nature of these (classical)
connections.
Finally, the covariant exterior derivative defined by the classical
connections
is a local differential operator in a very precise sense. If two fieds agree
on a a neighborhood, so do their covariant deriviatives. The connection we are
considering here is nonlocal also in that the basic "differential" operators
associated with it are not local in the sense just mentioned.

Given a nonlocal connection $\omega,$ we define a \textbf{derivative $D\psi$
with
respect to $\omega$} of a $p$-form $\psi,$ and then define the curvature of
$\omega$ to
be $D\omega.$ A nonlocal connection $\omega$ has a local, classical, aspect.
It is a usual, infinitesimal, connection and we denote it by the same letter
$\omega.$ We call a nonlocal connection \textbf{flat} if its local aspect is
flat.
Then we have a theorem which says that a nonlocal connection $\omega$ is flat
if and only if $D\omega = 0.$ Also, a nonlocal analog of Bianchi identity
holds:
$D^2\omega = 0.$

We also have a notion of plain \textbf{derivative} of a $p$-form without any
reference to a connection. This derivative is denoted by $d.$ It satisfies
$d^2 = 0.$ A nonlocal notion of integration is introduced, and a nonlocal analog of the
Stoke's theorem is also obtained. The latter is perhaps best described as
\textbf{the fundamental theorem of nonlocal calculus.}

All these notions can be made precise in the language of bundles. In
particular,
we can define \textbf{$p$-forms with values in an arbitrary vector bundle}.
Also, all the
notions of differential geometry have obvious analogs in this nonlocal
setting.
Thus, besides connections, we also have nonlocal metric structures,
$G$-structures, nonlocal holonomy, nonlocal structure equations, etc.

\section{Analysis} \label{analysis}
Analysis in the nonlocal setting takes two forms.

The first analyzes the usual local (point-wise) fields using nonlocal
connections. This entails the introduction of a nonlocal calculus of local
fields. (We call this \textbf{semilocal analysis:} local fields, nonlocal 
analysis.)
Here, too, we have natural notions of (nonlocal) integral and derivative. These
two are related by \textbf{the fundamental theorem of semilocal calculus,}
which says that these
two notions are inverses of each other. Like its Newtonian predecessor, it is
\emph{the} core of the calculus. Then, we analyze nonlocal differential
operators, on spaces of sections of vector bundles. The theory is very similar
in spirit to the usual analysis. Also,
the linear semilocal analysis is closely linked to the \emph{nonlinear} analysis
of the classical, local kind. This is, perhaps, one of the most promising
directions for future investigations.

The second one analyzes the nonlocal fields (i.e. $p$-forms), and as such is
completely nonlocal: nonlocal fields, nonlocal analysis. Here, the concepts,
and the machinery, are much different than those of local analysis. For one
thing, the algebra of $p$-forms is a graded \emph{noncommutative} ring. (This
fact
is
fruitfully exploited in the third part of \cite{nang}, described in the next
section.) Among other things, we study linear operators, and a detailed theory
of \textbf{elliptic} operators is obtained.

\section{Geometro-analysis} \label{synthesis}
The theories and concepts mentioned in the preceding sections are put together
to build nonlocal-noncommutative analogs of several classical theories such as
homology, cohomology, Hodge-de Rham, Chern-Weil, $K$-theory, index theory and
such (\cite{nang}, Part III). We re-emphasize that all these theories involve
functors with values in the category of \emph{noncommutative} rings. Naturally,
these theories encode much more information than do their classical,
commutative, counterparts. For each one of these theories, there is a natural
transformation from the noncommutative to the commutative case. Furthermore,
these natural transformations form commutative diagrams along with other
internal natural transformations of each theory. For example, the Chern
characters of nonlocal and local types form a commutative square of natural
transformations along with those relating the respective $K$-functors and the
$H$-functors, i.e. the cohomology theories.
This network of natural transformations
allow lifting problems---and often entire solutions---from the commutative to
the noncommutative functors on one hand while allowing translation of
noncommutative problems into their commutative analogs on the other. Of course,
this latter is more convenient in answering questions which have negative
answers.

\section{Algebraic Geometry} \label{ag}
The nonlocal methods can be easily extended to algebraic geometry and analytic
geometry: we have nonlocal sections of schemes giving noncommutative rings.
Thus, it becomes possible to analyze the geometry of noncommutative
rings/schemes using commutative schemes and sheaves. On the other hand,
nonlocal methods can shed
more light on commutative schemes. Thus, we have noncommutative Chow ring of a
scheme, noncommutative Grothendieck ring, nonlocal Grothendieck-Riemann-Roch,
etc.(\cite{nang}, Part IV). Here, too, we have a set of natural transformations
from the noncommutative functors to their classical commutative namesakes.

\section{Several Complex Variables} \label{ang}
The remarks of last section also apply to the category of analytic spaces and
the associated sheaves of modules. Nonlocal cohomology
of sheaves provide a much wider context in which to explore classical
theories, such as Oka-Cartan theory and their extensions/applications. Again,
we have a natural transformation from the 
noncommutative to the commutative cohomology (\cite{nang}, Part V).

\section{Number Theory} \label{nt}
This is, perhaps, the most difficult part of mathematics, and the nonlocal
viewpoint promises to be quite an addition to its already spectacular arsenal.
A lot of the arithmetic geometry has nonlocal extensions, not unlike the cases
of algebraic and analytic geometries discussed above. This is the field the
present author can not even pretend to be any proficient in; this is good news,
because the hope is that other truly qualified practitioners of this field
will find many more uses of the nonlocal viewpoint. Nevertheless, we mention
nonlocal Arakelov theory, and nonlocal arithmetic Riemann-Roch.

\end{document}